\documentclass[preprint,12pt]{article}

\usepackage{amssymb}
\usepackage{fullpage}
\usepackage{caption}
\usepackage{amsthm}
\usepackage{amsmath}
\usepackage{amssymb}
\usepackage{wasysym}
\usepackage{graphicx}
\usepackage{subcaption}
\usepackage{ragged2e}
\usepackage{float}
\usepackage{tikz}
\usetikzlibrary{decorations.pathmorphing,shadings, automata,positioning, decorations.markings, fixedpointarithmetic, shapes, decorations.pathreplacing, arrows, calc,math, arrows.meta,backgrounds}

\usepackage{verbatim}  % This includes a \begin{comment} environment
\usepackage{tabularx}
\usepackage{array,multirow,multicol}
\usepackage{subcaption}
\usepackage[super]{nth}
\usepackage{setspace}
\usepackage[ruled,vlined]{algorithm2e}
\usepackage{ algorithmicx}
\usepackage[a-1b]{pdfx}
\usepackage{hyperref}
\usepackage{xmpincl}

\tikzset{normal/.style={circle,draw, minimum size=8pt, inner sep=0,fill=black}}

\newtheoremstyle{strat}% name of the style to be used
  {\topsep}% measure of space to leave above the theorem. E.g.: 3pt
  {\topsep}% measure of space to leave below the theorem. E.g.: 3pt
  {}% name of font to use in the body of the theorem
  {0pt}% measure of space to indent
  {\bfseries}% name of head font
  {.}% punctuation between head and body
  { }% space after theorem head; " " = normal interword space
  {\thmname{#1}\thmnumber{ #2}\textnormal{\thmnote{ (#3)}}}

\newtheorem{theorem}{Theorem}

\numberwithin{equation}{section}

\newtheorem{lemma}[theorem]{Lemma}
\newtheorem{corollary}[theorem]{Corollary}

\newtheorem{conj}[theorem]{Conjecture}
{\theoremstyle{strat}
\newtheorem{strat}[theorem]{Strategy}}
\theoremstyle{definition}

\DeclareUnicodeCharacter{2212}{-}
\numberwithin{theorem}{section}

\title{Brushing Directed Graphs}

\author{Jared Howell\thanks{School of Science and the Environment, Grenfell Campus, Memorial University of Newfoundland, Corner Brook, NL, A2H 5G5, Canada.}
\and
Sulani D.~Kavirathne\thanks{Department of Mathematics and Statistics, Memorial University of Newfoundland, St.~John's, NL, A1C 5S7, Canada.}
\and
David A.~Pike\thanks{Department of Mathematics and Statistics, Memorial University of Newfoundland, St.~John's, NL, A1C 5S7, Canada. \texttt{dapike@mun.ca}}
}

\begin{document}
\maketitle

\begin{abstract}
Brushing of graphs is a graph searching process in which the searching agents are called brushes.  We focus on brushing directed graphs based on a new model in which the brushes can only travel in the same direction as the orientation of the arcs that they traverse.
%The fundamentals of this model have been influenced by previous studies about brushing undirected graphs.
We discuss strategies to brush directed graphs as well as values and bounds for the brushing number of directed graphs. We determine the brushing number for any transitive tournament, which we use to give an upper bound for the brushing number of directed acyclic graphs in general.  We also establish exact values for the brushing numbers of complete directed graphs, rooted trees, and rotational tournaments.
%The behaviour of the brushing number of several other types of directed graphs are also taken into consideration.
\end{abstract}



\iffalse
\fi


\section{Introduction}
\label{sec1}
As noted in a survey by Alspach~\cite{alspach2004searching},
there are several types of graph searching that can be used to search networks, depending on particular search objectives and paradigms.  Examples include such goals as capturing a visible intruder (as is modelled by the game of Cops and Robber~\cite{BonatoNowakowski2011}), or
conducting periodic patrols to detect nearby intruders (such as is accomplished by a watchman's walk~\cite{DHP2021,HRW1998}).
These models often make use of some resource that is sought to be optimised, typically either personnel or time, that must be allocated to the search effort.
In this paper we focus on a type of graph searching that can be used to locate an invisible intruder (if one is present) in a network by means of a coordinated search of the entire network.  The model that we consider is itself a type of graph cleaning; it assumes that a given network has been infiltrated, and deploys a collection of cleaning agents that travel throughout the corresponding graph in a manner that ultimately assures that all parts of the graph are guaranteed to be free of intrusion.
Specifically, we consider the graph searching technique known as brushing, and adapt it to the setting of directed graphs.

Graph brushing is a type of graph searching that draws inspiration from chip firing.
In chip firing~\cite{BJORNER1991283} there is an initial configuration of tokens called \textit{chips} that are placed on vertices.  A vertex is said to be \textit{primed} if the vertex has at least as many chips as its degree. When a primed vertex \textit{fires}, it sends a chip along each incident
edge, with the neighbouring vertices accumulating the dispersed chips.
%In the study of chip firing, commonly raised questions have been variants of “is this process finite or infinite?”, “how many chips are needed to produce a cycle?” and “how long will it take for an infinite game to become a repeated cycle of chip configurations?”.
In the context of brushing, the tokens are called \textit{brushes}, and we treat the entire graph as being in an initial state of dirtiness.
As an analogy consider a network of pipes and junctions that have to be cleaned of a contaminant such as algae.
When a vertex fires its brushes, the vertex becomes clean, as do the incident edges that the brushes travel along.
%Instead of requiring a vertex $v$ to have at least $\deg(v)$ brushes before it can fire, it now suffices for $v$ to have at least as many brushes as incident edges that have not yet been cleaned.

The process of brushing undirected graphs was introduced by McKeil~\cite{McKeil2007}.
%When a vertex becomes \textit{cleaned}, a brush must travel down each incident contaminated edge. Once a brush traverses an edge, that edge is \emph{cleaned}.
The ultimate goal is to clean a given graph $G$, which is accomplished once every edge and every vertex of $G$ has been cleaned.  The location of the brushes will change over time, so at each discrete time step there is a \textit{brush configuration} that records the positions of the brushes and the state of cleanliness of each vertex and edge.  McKeil discussed two distinct objectives, one being how many brushes are needed to clean a graph,
and the other being the question of how quickly a graph can be cleaned (under the premise that multiple vertices can disperse their brushes in parallel).
%In \cite{McKeil2007}  \textit{dispersal mode} was defined as the rule that chooses the way that vertices can fire.  Two types of dispersal modes, namely parallel and sequential, are given. McKeil developed two models, called the regular cleaning model and open cleaning model.  For both models McKeil presented strategies to minimise the number of brushes used and number of time steps taken to clean a graph. Also in the two models more than one brush can travel through an edge and brushes can travel through cleaned edges.

Brushing undirected graphs has received subsequent attention and study, for which the \textit{brushing number} has been the primary focus of interest, namely the least number of brushes required by any initial brush configuration that enables a given graph to be successfully cleaned.
During the cleaning process, instead of requiring a vertex $v$ to have at least $\deg(v)$ brushes before it can fire its brushes, it now suffices for $v$ to have at least as many brushes as incident edges that have not yet been cleaned.

A model in which brushes cannot traverse clean edges and each dirty edge must be traversed by only one brush during the cleaning process
is considered by Messinger, Nowakowski and Prałat~\cite{messinger2008cleaning}.
%whereby the number of brushes that can traverse each edge is limited.
\iffalse
In that model every edge in a graph $G$ is initially considered dirty and a fixed number of brushes begin on a set of
vertices.
At each step of the process, a vertex $v$ may be cleaned
%(instead of fired)
if the number of brushes on $v$ is greater than or equal to the number of dirty incident edges.
\fi
A less restrictive model is considered by Bryant, Francetić, Gordinowicz, Pike and Prałat in~\cite{bryant2014brushing}, where dirty edges are allowed to be traversed by multiple brushes.  With respect to this model, general bounds are established for the brushing number of an undirected graph, including bounds that are given in terms of parameters
such as cutwidth and bisection width.
Trees are also considered in~\cite{bryant2014brushing}, whereby it is proved that
for a tree $T$ with $\ell(T)$ leaves, $\frac{\ell(T) + 1}{2}$ is the brushing number when $\ell(T)$  is odd.  However, when $\ell(T)$ is even the
brushing number is either $\frac{\ell(T)}{2}$ or $\frac{\ell(T)}{2} + 1$, implying that trees with an even number of leaves are divided into two groups depending on their brushing number.   Penso, Rautenbach and Ribeiro de Almeida subsequently proved that trees with  $\ell(T)$ leaves and brushing number $\frac{\ell(T)}{2}$ can be recognised efficiently~\cite{MR3391773Penso}.

In the present paper (which is based on thesis given in~\cite{Sulani2024}) we generalise the notion of brushing to directed graphs.
In Section~\ref{sec p2}, we formulate the basic rules of brushing for directed graphs, along with defining necessary concepts and notation. In Section~\ref{sec 2}, the brushing number of directed acyclic graphs is discussed in detail.
We prove that every transitive tournament on $n \geq 3$ vertices has brushing number $\lfloor \frac{n^2}{4} \rfloor$ (see Theorem~\ref{th_1.3}).  We use this result to establish an upper bound on the brushing number of directed acyclic graphs in general (see Theorem~\ref{th:6}). In Section~\ref{chap 3}, we take several other types of directed graphs into consideration and exactly determine the brushing numbers of complete directed graphs, rooted trees, and rotational tournaments. We also show that other concepts like path decomposition and the transpose of a directed graph are associated with the brushing number. In the final section we present some discussion and questions for future inquiry.

\section{Some Preliminaries}
\label{sec p2}
We adapt the model in \cite{bryant2014brushing} by making some modifications relevant for directed graphs.
In the process of cleaning a directed graph $G$, an isolated vertex can fire provided that it has a brush, whereas a vertex $v\in V(G)$ that is not isolated is able to fire only if the present number of brushes at $v$ is at least the out-degree of $v$. When a vertex $v$ fires, the brushes on $v$ clean $v$, and each outgoing arc incident with $v$ is traversed by at least one brush that is fired from $v$, thereby cleaning the dirty outgoing arcs that were incident with $v$. At the end of the time step, each brush from the firing vertex $v$ moves to the vertex adjacent to $v$ at the endpoint of the arc it traversed. Excess brushes are allowed to remain at $v$. Arcs are allowed to be traversed by more than one brush, and multiple brushes can simultaneously traverse an arc; however each arc can be traversed during at most one step of the cleaning process. We allow vertices to fire one at a time in a sequence. In this model we do not fire multiple vertices that are ready to fire at the same time step (although a way of simultaneously firing vertices when brushing undirected graphs is described as \textit{parallel dispersal mode}  in \cite{McKeil2007}). Before beginning the process of cleaning we decide the sequence in which the vertices fire. Since we allow multiple brushes to traverse an arc, it is possible that a vertex which fires later in the sequence might accumulate brushes from vertices which have been fired earlier in the sequence.

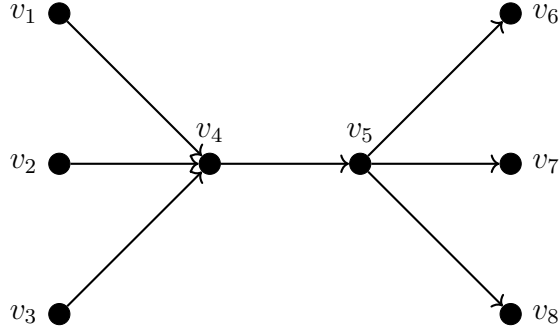
\begin{figure}[h]
    \centering
    \begin{tikzpicture}

         \node[normal] (v1) at (1,-2)[label=left:$v_1$]{};
        \node[normal] (v2) at (1,-4) [label=left:$v_2$]{};
        \node[normal] (v3) at (1,-6)[label=left:$v_3$] {};
        \node[normal] (v4) at (3,-4)[label=above:$v_4$] {};
        \node[normal] (v5) at (5,-4)[label=above:$v_5$] {};
         \node[normal] (v6) at (7,-2) [label=right:$v_6$]{};
         \node[normal] (v7) at (7,-4) [label=right:$v_7$]{};
         \node[normal] (v8) at (7,-6) [label=right:$v_8$]{};

       \draw[thick,->] (v1) to (v4);
       \draw[thick,->] (v2) to (v4);
       \draw[thick,->] (v3) to (v4);
       \draw[thick,->] (v4) to (v5);

       \draw[thick,->] (v5) to (v6);
      \draw[thick,->] (v5) to (v7);
      \draw[thick,->] (v5) to (v8);

    \end{tikzpicture}
    \captionsetup{justification=raggedright, margin=0.8cm}
    \caption{An example of a directed graph where multiple brushes simultaneously traverse an arc during brushing.}
    \label{fig:11}
\end{figure}

As an example refer to the directed graph $G$ in Figure~\ref{fig:11}. At the beginning of the process of cleaning this graph $G$, place one brush at each of the vertices $v_1,v_2$ and $v_3$. Then $v_1,v_2$ and $v_3$ fire sequentially and $v_4$ accumulates three brushes. Next $v_4$ fires and these three brushes traverse the arc $(v_4,v_5)$. Then $v_5$ receives three brushes and $v_5$ fires, resulting in the vertices $v_6,v_7$ and $v_8$ getting one brush each. Then $v_6,v_7$ and $v_8$ fire sequentially; brushes stay at these vertices as these are sink vertices. Thus we have cleaned the graph $G$ with three brushes.

In general, the process of cleaning terminates when a clean graph is obtained, or else the process stops in a situation in which there are dirty vertices but none are capable of firing. Let $B^t_{G}(v)$ denote the number of brushes at vertex $v$ of $G$, at time $t$ ($t=0,1,2, \ldots, |V(G)|$), so that $B^0_{G}(v)$ denotes the number of brushes at vertex $v$ in the initial configuration. A vertex $v_i$ fires in between time $i-1$ and $i$. The \emph{brushing number} of a graph is the minimum number of brushes needed for some initial configuration to clean the graph. For a graph $G$, the brushing number is denoted by $B(G)$.

\begin{figure}[h]
    \centering
\begin{tikzpicture}
%\tikzstyle{every node}=[draw, shape=circle];
%one
\node[circle] (v0) at (0,0) [fill,label=below:$v_2$] {};
\node[circle,inner sep=2pt,draw,text=white,font=\tiny] (v1) at (1,0) [fill,label=below:$v_1$] {\textbf{2}};
\node[circle] (v2) at (0,2) [fill,label=above:$v_3$] {};
\node[circle] (v3) at (1,2) [fill,label=above:$v_4$] {};
\draw[thick] (v1) to (v0) node[ above right = -1.8cm of v1] {$t=0$};
\draw[thick] (v1) to (v2);
\draw[thick] (v0) to (v2);
\draw[thick] (v2) to (v3);
%two
\node[circle,draw] (v1) at  (4,0) [label=below:$v_1$] {};
\node[inner sep=2pt,draw,font=\tiny,shape=circle,fill,text=white](v0) at (3,0) [fill,label=below:$v_2$] {\textbf{1}};
\node[inner sep=2pt,draw,font=\tiny,shape=circle,fill,text=white] (v2) [fill,label=above:$v_3$] at (3,2){\textbf{1}};
\node[circle] (v3) at (4,2) [fill,label=above:$v_4$] {};
\draw[dashed,thick] (v1) to (v0) node[ above right = -1.8cm of v1] {$t=1$};
\draw[dashed,thick] (v1) to (v2);
\draw[thick] (v0) to (v2);
\draw[thick] (v2) to (v3);
%three
\node[circle,draw] (v1) at  (7,0) [label=below:$v_1$] {};
\node[circle,draw]  (v0) at  (6,0) [label=below:$v_2$] {};
\node[inner sep=2pt,draw,font=\tiny,shape=circle,fill,text=white](v2) [label=above:$v_3$] at (6,2) {\textbf{2}};
\node[circle] (v3) at (7,2) [fill,label=above:$v_4$] {};
\draw[dashed,thick] (v1) to (v2);
\draw[dashed,thick] (v1) to (v0) node[ above right = -1.8cm of v1] {$t=2$};
\draw[dashed,thick] (v0) to (v2);
\draw[thick] (v2) to (v3);
%four

\node[circle,draw]  (v1) at  (10,0) [label=below:$v_1$] {};
\node[circle,draw]  (v0) at  (9,0) [label=below:$v_2$] {};
\node[text=black,inner sep=2pt,draw,font=\tiny,circle] (v2) at (9,2)[label=above:$v_3$] {\textbf{1}};
\node[inner sep=2pt,draw,font=\tiny,shape=circle,fill,text=white] (v3) at (10,2) [label=above:$v_4$] {\textbf{1}};
\draw[dashed,thick] (v1) to (v2);
\draw[dashed,thick] (v1) to (v0) node[ above right = -1.8cm of v1] {$t=3$};
\draw[dashed,thick] (v0) to (v2);
\draw[dashed,thick,] (v2) to (v3);
%five
\node[circle,draw]  (v1) at  (13,0) [label=below:$v_1$] {};
\node[circle,draw]  (v0) at  (12,0) [label=below:$v_2$] {};
\node[text=black,inner sep=2pt,draw,font=\tiny,circle] (v2) at (12,2)[label=above:$v_3$] {\textbf{1}};
\node[text=black,inner sep=2pt,draw,font=\tiny,circle] (v3) at (13,2) [label=above:$v_4$] {\textbf{1}};
\draw[dashed,thick] (v1) to (v2);
\draw[dashed,thick] (v1) to (v0) node[ above right = -1.8cm of v1] {$t=4$};
\draw[dashed,thick] (v0) to (v2);
\draw[dashed,thick] (v2) to (v3);
\end{tikzpicture}
 \caption{Graph $G^{'}$ with brush number 2.}
    \label{fig:12}
\end{figure}
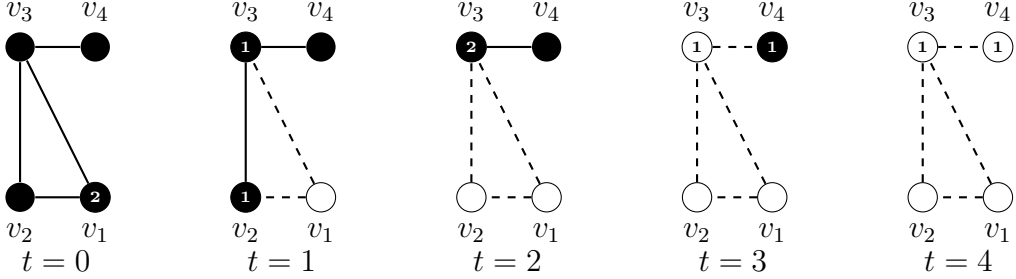

As an example refer to Figure~\ref{fig:9}. The process of brushing the underlying undirected graph $G^{'}$ of the directed graph $G$ in Figure~\ref{fig:9} with the cleaning model in \cite{messinger2008cleaning} is given in Figure~\ref{fig:12}. In Figures~\ref{fig:9} and \ref{fig:12}, black nodes represent unfired (or dirty) vertices. The edges and arcs in thick lines represent dirty edges and dirty arcs. The white nodes represent fired (or clean) vertices. The edges and arcs in dashed lines represent clean edges and clean arcs. The number of brushes at a vertex in a time step is given inside the node. If there is no number inside a vertex then that vertex does not have any brushes at the particular time step.

Based on the cleaning model that we have described for a directed graph $G$ with arc set $A(G)$, the following theorem establishes initial bounds on the brushing number $B(G)$.
\begin{theorem}\label{th 1}
    If $G$ is a directed graph, then $\max\{\deg ^+(v)\mid v\in V(G)\}\le B(G)\le |A(G)|.$
\end{theorem}
\begin{proof}
    In a directed graph $G$ if $v_M$ is a vertex with maximum out-degree, to fire $v_M$ we need $\deg^+(v_M)$ brushes. Therefore $\deg^+(v_M)\le B(G)$. Also, when all the vertices have been fired, then at least one brush should have traversed each arc, which implies $|A(G)|$ is an upper bound for the brushing number.\end{proof}
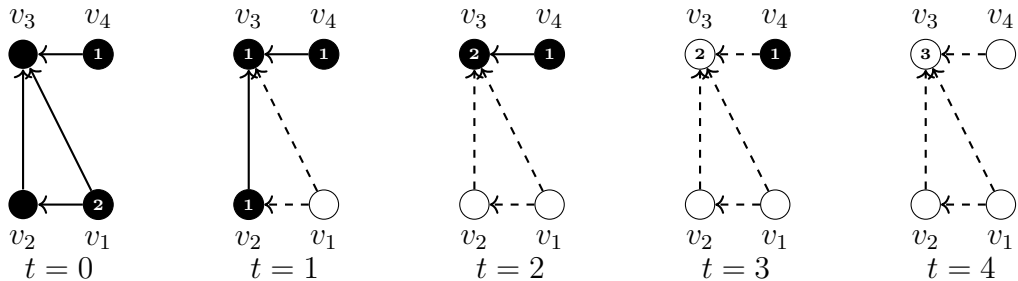
\begin{figure}[h]
    \centering
\begin{tikzpicture}
%\tikzstyle{every node}=[draw, shape=circle];
%one
\node[circle] (v0) at (0,0) [fill,label=below:$v_2$] {};
\node[circle,inner sep=2pt,draw,text=white,font=\tiny] (v1) at (1,0) [fill,label=below:$v_1$] {\textbf{2}};
\node[circle] (v2) at (0,2) [fill,label=above:$v_3$] {};
\node[circle,inner sep=2pt,draw,text=white,font=\tiny] (v3) at (1,2) [fill,label=above:$v_4$] {\textbf{1}};
\draw[thick,->] (v1) to (v0) node[ above right = -1.8cm of v1] {$t=0$};
\draw[thick,->] (v1) to (v2);
\draw[thick,->] (v0) to (v2);
\draw[thick,->] (v3) to (v2);
%two
\node[circle,draw] (v1) at  (4,0) [label=below:$v_1$] {};
\node[inner sep=2pt,draw,font=\tiny,shape=circle,fill,text=white](v0) at (3,0) [fill,label=below:$v_2$] {\textbf{1}};
\node[inner sep=2pt,draw,font=\tiny,shape=circle,fill,text=white] (v2) [fill,label=above:$v_3$] at (3,2){\textbf{1}};
\node[circle,inner sep=2pt,draw,text=white,font=\tiny] (v3) at (4,2) [fill,label=above:$v_4$] {\textbf{1}};
\draw[dashed,thick,->] (v1) to (v0) node[ above right = -1.8cm of v1] {$t=1$};
\draw[dashed,thick,->] (v1) to (v2);
\draw[thick,->] (v0) to (v2);
\draw[thick,->] (v3) to (v2);
%three
\node[circle,draw] (v1) at  (7,0) [label=below:$v_1$] {};
\node[circle,draw]  (v0) at  (6,0) [label=below:$v_2$] {};
\node[inner sep=2pt,draw,font=\tiny,shape=circle,fill,text=white](v2) [label=above:$v_3$] at (6,2) {\textbf{2}};
\node[circle,fill,inner sep=2pt,draw,text=white,font=\tiny](v3) at (7,2) [label=above:$v_4$]{\textbf{1}};
\draw[dashed,thick,->] (v1) to (v2);
\draw[dashed,thick,->] (v1) to (v0) node[ above right = -1.8cm of v1] {$t=2$};
\draw[dashed,thick,->] (v0) to (v2);
\draw[thick,->] (v3) to (v2);
%four

\node[circle,draw]  (v1) at  (10,0) [label=below:$v_1$] {};
\node[circle,draw]  (v0) at  (9,0) [label=below:$v_2$] {};
\node[text=black,inner sep=2pt,draw,font=\tiny,circle] (v2) at (9,2)[label=above:$v_3$] {\textbf{2}};
\node[inner sep=2pt,draw,font=\tiny,shape=circle,fill,text=white] (v3) at (10,2) [label=above:$v_4$] {\textbf{1}};
\draw[dashed,thick,->] (v1) to (v2);
\draw[dashed,thick,->] (v1) to (v0) node[ above right = -1.8cm of v1] {$t=3$};
\draw[dashed,thick,->] (v0) to (v2);
\draw[dashed,thick,->] (v3) to (v2);
%five
\node[circle,draw]  (v1) at  (13,0) [label=below:$v_1$] {};
\node[circle,draw]  (v0) at  (12,0) [label=below:$v_2$] {};
\node[text=black,inner sep=2pt,draw,font=\tiny,circle] (v2) at (12,2)[label=above:$v_3$] {\textbf{3}};
\node[draw,circle] (v3) at (13,2) [label=above:$v_4$] {};
\draw[dashed,thick,->] (v1) to (v2);
\draw[dashed,thick,->] (v1) to (v0) node[ above right = -1.8cm of v1] {$t=4$};
\draw[dashed,thick,->] (v0) to (v2);
\draw[dashed,thick,->] (v3) to (v2);
\end{tikzpicture}
 \caption{A directed graph $G$ with $B(G)= 3$.}
    \label{fig:9}
\end{figure}

\section{Brushing Directed Acyclic Graphs }
\label{sec 2}
This section is about brushing transitive tournaments and directed acyclic graphs. A \textit{tournament} is a directed graph in which there is exactly one arc between each pair of vertices. A tournament  $\mathcal{T} $ is \textit{transitive} if and only if for all vertices $u, v, w\in V(\mathcal{T})$ if $(u,v)$ and $(v,w)$ are arcs of $\mathcal{T}$ then $(u,w)$ is an arc of $\mathcal{T}$. Note that if $(u,v)$ and $(v,w)$ are arcs in a tournament $\mathcal{T}$, then we must have $w\neq u$, and $(w,u)$ is an arc of $\mathcal{T}$ if and only if $(u,w)$ is not an arc of $\mathcal{T}$.

\subsection{Transitive tournaments}\label{sec1}
In this subsection we begin by introducing a strategy that can be used to brush a transitive tournament, followed by a lemma that will be used to find the brushing number of a transitive tournament.  Further, the  brushing number of a transitive tournament when one edge is taken away will be discussed.

\begin{strat}\label{strat1}
  \begin{justify}
In this strategy, a series of steps for brushing a transitive tournament with $n \geq 3$ vertices will be presented.
Consider the set of vertices $\{v_1,v_2, v_3,\dots, v_n\}$ with  $\deg^+(v_i)=n-i,\:\mbox{for all}\;1\le i\le n$. Let $B^t_{G}(v)$ be the number of brushes at vertex $v$ of $G$, at time $t$ ($t=0,1,2, \ldots, n$). For each $i\in {1,2,\ldots,n}$ the vertex $v_i$ fires in between time $i-1$ and $i$. $\mbox{Let}\:\ell =\lfloor\frac{n}{2}\rfloor$. Let $B^0_{G}(v_{k})$ denote the number of brushes at the $k^{\text{th}}$ vertex in the initial configuration such that
$$B^0_{G}(v_{k}) =
    \begin{cases}
     n-(2k-1) & \text{if }  k=1,2,\ldots, \ell,\\
        0 & \text{}  \mbox{otherwise.}
    \end{cases}$$
    So no vertex $v_i$ with index $i$ greater than $\ell$ receives a brush in the initial configuration. As vertex $v_1$ with out-degree $n-1$ has  $n-1$  brushes,  $v_1$ fires.  Now

$$B^1_{G}(v_{k}) =
    \begin{cases}
    0 & \text{if } k=1, \\
     n-(2k-2) & \text{if }  k=2, 3,\ldots, \ell,\\
        1 & \text{}  \mbox{otherwise.}
    \end{cases}$$
Observe that the number of brushes at $v_2$ is now $n-2$. So $v_2$ fires.  Then
$$B^2_{G}(v_{k}) =
    \begin{cases}
    0 & \text{if }  k=1,2,\\
     n-(2k-3
     ) & \text{if }  k=3,4,\ldots, \ell,\\
        2 & \text{}  \mbox{otherwise.}
    \end{cases}$$
After the $i^{\text{th}}$ vertex fires $(i=1,2,\ldots,\ell)$
$$B^i_{G}(v_{k}) =
    \begin{cases}
    0
      & \text{if }  k=1, 2,\ldots, i,\\
     n-2k+i+1
      & \text{if }  k=i+1,i+2,\ldots, \ell,\\
        i & \text{}  \mbox{otherwise.}
    \end{cases}$$
When $v_1, v_2,\ldots,v_{i-1}$ have all fired $v_i$ is able to fire, as it will then have $n-(2i-1)+i-1=n-i=\deg^+(v_i)$ brushes, where $i=1,2,\ldots,\ell$.
After the vertices $v_1,v_2,\ldots,v_{\ell-1}$ fire, $v_\ell$ obtains $\ell-1$ brushes. But initially, there were $n-(2\ell-1)$ brushes in $v_\ell$. Therefore, in total $v_\ell$ has $n-\ell$ brushes, which is equal to its out-degree. Then $v_\ell$ fires. Initially there were no brushes in $v_{\ell+1}$. After $v_\ell$ fires, $v_{\ell+1}$ has accumulated $\ell$ brushes. As $\deg^+(v_{\ell+1})=n-(\ell+1)$, $\deg^+(v_{\ell+1})= \lfloor\frac{n-1}{2}\rfloor$ and then $n-(\ell+1) = \deg^+(v_{\ell+1})\le \ell$. Therefore $v_{\ell+1}$ fires. Similarly, we can  also show that $v_{\ell+2},v_{\ell+3},\ldots, v_n$ are able to fire in sequence. This concludes Strategy~\ref{strat1}.
\end{justify}
\end{strat}

In order to prove the equality in Theorem~\ref{th_1.3} we will need Lemma~\ref{le_1.5}, for which  we will introduce some notations. Let $G$ be a directed graph. The set of edges with one end in $S$ and the other end in $T$, where $S,T \subseteq V(G)$, is denoted by $[S,T]$. An \textit{edge-cut} is a set of edges of the form $[S,\Bar{S}]$ where $\emptyset \subset S\subset V(G)$ and $\Bar{S}=\{v\in V(G) \mid v\notin S\}$. The directed subgraph of $G$ whose vertex set is $S$ and whose arc set is the set of those arcs of $G$ that have both ends in $S$ is called the \textit{subgraph of $G$ induced by $S$} and is denoted by $G[S]$.

\begin{lemma}\label{le_1.5}
If $G$ is a directed graph and $S\subseteq V(G)$ such that no arcs go from $\Bar{S}$ to $S$, then $B(G)\ge|[S,\Bar{S}]|$. Moreover, when $G$ is a transitive tournament $B(G)\ge|S||\Bar{S}|$.
\end{lemma}
\begin{proof}
    Let $G$ be a directed graph and $S\subseteq V(G)$ such that no arcs go from $\Bar{S}$ to $S$. As a result of the edge-cut $[S,\Bar{S}]$ two subgraphs of $G$,  $G[S]$ and $G[\Bar{S}]$, are obtained. The subgraphs $G[S]$ and $G[\Bar{S}]$ can be brushed as two separate graphs. But  consider $G$ where $A(G)=A(G[S])\cup A(G[\Bar{S}])\cup [S,\Bar{S}]$ and the three sets of arcs $A(G[S]),\,A(G[\Bar{S}])$ and $[S,\Bar{S}]$ are pairwise disjoint. When brushing $G$, we have to use one brush for each arc in $[S,\Bar{S}]$. Therefore $B(G)\ge|[S,\Bar{S}]|$.

    When $G$ is a transitive tournament, consider a set $S\subseteq V(G)$ such that no arcs go from $\Bar{S}$ to $S$. Then $|[S,\Bar{S}]|=|S||\Bar{S}|$ as every vertex in $S$ is adjacent to every vertex in $\Bar{S}$. Therefore $B(G)\ge|S||\Bar{S}|$.
\end{proof}
\begin{theorem}\label{th_1.3}
If $G$ is a transitive tournament with $n \geq 3$ vertices, then
$$B(G)=
    \begin{cases}
     \frac{n^2-1}{4} & \text{if } \mbox{n is odd,}\\
        \frac{n^2}{4} & \text{if } \mbox{n is even.}
    \end{cases}$$
\end{theorem}

\begin{proof}
When Strategy~\ref{strat1} is used to brush $G$,
$ B(G)\le\sum\limits_{k=1}^{\ell} B^0_{G}(v_{k})$.
So when $n$ is even, $$B(G) \le\underbrace{ n-1+n-3+\cdots+1}_\text{sum of odd numbers between $0$ and $n$} = \left(\frac{n}{2}\right)^2 = \frac{n^2}{4},$$
and when $n$ is odd, $$B(G) \le\underbrace{ n-1+n-3+\cdots+2}_\text{sum of even numbers between $0$ and $n$}
= \left(\frac{n-1}{2}\right)\left(\frac{n-1}{2}+1\right)=\frac{(n^2-1)}{4}.$$

By using Lemma~\ref{le_1.5} for a transitive tournament $G$ with $S=\{v_1,v_2, v_3,\dots, v_\ell\}$, we get the following lower bound on $B(G)$:
$$B(G)\ge
    \begin{cases}
     \frac{n^2-1}{4} & \text{if } \mbox{$n$ is odd,}\\
        \frac{n^2}{4} & \text{if } \mbox{$n$ is even.}
    \end{cases}$$
    Therefore
$$B(G)=
    \begin{cases}
     \frac{n^2-1}{4} & \text{if } \mbox{$n$ is odd,}\\
        \frac{n^2}{4} & \text{if } \mbox{$n$ is even.}
    \end{cases}$$
\end{proof}

An immediate consequence of Theorem~\ref{th_1.3} is that
%The obtained brushing number for a transitive tournament implies that
Strategy~\ref{strat1} is an optimal strategy when applied to transitive tournaments.
By observing the movement of brushes in a transitive tournament as described by this strategy, we obtain the following corollary.

\begin{corollary}\label{co_1}
If $G$ is a transitive tournament that is cleaned using Strategy~\ref{strat1}, then there is at least one brush that travels a Hamiltonian path during the brushing of $G$.
\end{corollary}
\begin{proof}
When using Strategy~\ref{strat1} to brush a transitive tournament $G$, every vertex obtains one brush from the preceding vertex that just fired.
\end{proof}

So far in this section we have discussed brushing transitive tournaments and facts related to this process. We now ask: if one of the arcs in a transitive tournament is removed, what will be the brushing number of the resulting directed graph?   We show that removing a single arc from a transitive tournament does not increase the brushing number.

\begin{theorem}\label{th_1.4}
Let $G$ be a transitive tournament on $n\ge 3$ vertices. If $e\in A(G)$, then $B(G-e)\le B(G)$.
\end{theorem}

\begin{figure}[h]
    \centering
\begin{tikzpicture}
%\tikzstyle{every node}=[draw, shape=circle];
%one
\node[circle,fill] (v0) at (0,0)  {};
\node[circle,fill] (v1) at (1.05,2) {};
\node[circle,fill] (v2) at (3,1.5)  {};
\node[circle,fill] (v3) at (5,1.5) {};
\node[circle,fill] (v4) at (6.95,2)  {};
\node[circle,fill] (v5) at (8,0) {};
\node[] (v6) at (6,1.8) {};
\draw[thick,->] (v0) to (v1);
\draw[thick,->] (v0) to (v2);
\draw[thick,->] (v1) to (v2);
\draw[thick,->] (v2) to (v3)node[left = 1.5cm of v6] {$e$};
\draw[thick,->] (v3) to (v4);
\draw[thick,->] (v4) to (v5);
\draw[thick,->] (v3) to (v5);
\draw[thick,->] (v0) to (v5);
\end{tikzpicture}
\caption{A directed graph $G$ such that $B(G) < B(G-e)$.}
    \label{fig:10}
\end{figure}
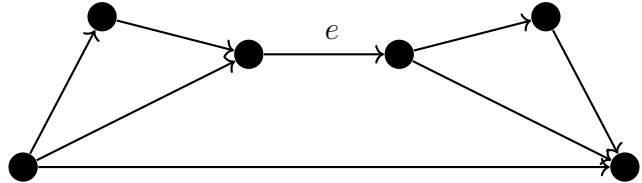

\begin{proof}

Label the vertices of $G$ in the descending order of their out-degree, so that
\\$$\deg^+_{G}(v_1)> \deg^+_G(v_2)> \cdots  > \deg^+_G(v_n).$$
\begin{comment}
  In the process of brushing, let $v_{L^\ast}$ be the last vertex we add brushes to. If n is odd, $v_{L^\ast}$ is the $\frac{(n-1)}{2}$ $^{th}$ vertex and if $n$ is even, $v_{L^\ast}$ is the $\frac{n}{2}$ $^{th}$ vertex, in the descending order of vertices.
\end{comment}

Let $\ell =\lfloor\frac{n}{2}\rfloor$. Strategy~\ref{strat1} will be used for cleaning  $G$ and based on that a strategy to clean $G-e$ using no more than $B(G)$ brushes will be given in this proof. Let $L=\{v_1, v_2, v_3, v_4,\ldots,v_{\ell}\}$  and
$R=V(G)\setminus L$.
Recall that $B^t_{G}(v)$ denotes the number of brushes at vertex $v$ of $G$, at time $t$ ($t=0,1,2, \ldots, n$) and $B^0_{G}(v_{i})$ denotes the number of brushes at the $i^{\text{th}}$ vertex in the initial configuration. The vertex $v_i$ fires in between time $i-1$ and $i$. Observe that for each $v_i \in R \setminus \{v_{\ell +1}\}$, $B^{i-1}_{G}(v_i) > \deg^{+}_{G}(v_i)$. We will describe an initial configuration $B^0_{G-e}$ to clean $G-e$, whereby the vertices will fire in the sequence $v_1, v_2, \ldots, v_n$.
There are several cases, depending on the nature of the vertices of the arc $e=(v_a,v_b)$ and the parity of $n$.

Case I: $v_a \in L \setminus \{v_{\ell}\}$ and $v_b \in R \setminus \{v_{\ell +1}\}$.
\\
As the arc $e$ is not present  in $G-e$,
let $B^0_{G-e}(v_a)=B^0_{G}(v_a)-1$  and  $B^0_{G-e}(v)=B^0_{G}(v)$ for each $v\in V(G) \setminus \{v_a\}$. Then  for $0\le t\le a-1$,
$$B^t_{G-e}(v) =
    \begin{cases}
     B^t_{G}(v)-1 & \text{if }  v=v_a,\\
        B^t_{G}(v) & \text{}  \mbox{otherwise.}
    \end{cases}$$
Now for $t\in\{a, a+1,\ldots, b-1\}$,
$$B^t_{G-e}(v) =
    \begin{cases}
     B^t_{G}(v)-1 & \text{if }  v=v_b,\\
        B^t_{G}(v) & \text{}  \mbox{otherwise.}
    \end{cases}$$
Observe that $B^{b-1}_{G-e}(v_b)=B^{b-1}_{G}(v_b)-1\ge \deg^{+}_{G}(v_b)=\deg^{+}_{G-e}(v_b)$  and so $v_b$ can fire between time $b$ and $b-1$. Consequently for $t\in\{b,b+1,\ldots, n\}$,
$$B^t_{G-e}(v) =
    \begin{cases}
     B^t_{G}(v)-1 & \text{if }  v=v_b,\\
        B^t_{G}(v) & \text{}  \mbox{otherwise.}
    \end{cases}$$
    Thus  $B(G-e)\le B(G)-1$.

Case II:  $v_a \in L \setminus \{v_{\ell}\}$, $v_b  = v_{\ell +1}$ and $n$ is even.
\\
As the arc $e$ is not present  in $G-e$,
let $B^0_{G-e}(v_a)=B^0_{G}(v_a)-1$  and  $B^0_{G-e}(v)=B^0_{G}(v)$ for each $v\in V(G) \setminus \{v_a\}$.
Then  for $t\in\{0,1,\ldots, a-1\}$,
$$B^t_{G-e}(v) =
    \begin{cases}
     B^t_{G}(v)-1 & \text{if }  v=v_a,\\
        B^t_{G}(v) & \text{}  \mbox{otherwise.}
    \end{cases}$$
Now for $t\in\{a, a+1,\ldots, \ell\}$,
$$B^t_{G-e}(v) =
    \begin{cases}
     B^t_{G}(v)-1 & \text{if }  v=v_{\ell +1},\\
        B^t_{G}(v) & \text{}  \mbox{otherwise.}
    \end{cases}$$
Observe that $B^{\ell}_{G-e}(v_{\ell+1})=B^{\ell}_{G}(v_{\ell+1})-1= \deg^+_{G-e}(v_{\ell+1})$  and so $v_{\ell+1}$ can fire. Consequently for $t\in\{\ell+1,\ell+2,\ldots, n\}$,
$$B^t_{G-e}(v) =
    \begin{cases}
     B^t_{G}(v)-1 & \text{if }  v=v_{\ell +1},\\
        B^t_{G}(v) & \text{}  \mbox{otherwise.}
    \end{cases}$$
Thus $B(G-e)= B(G)-1$.

Case III:  $v_a \in L \setminus \{v_{\ell}\}$, $v_b = v_{\ell +1}$ and $n$ is odd.
\\
As the arc $e$ is not present  in $G-e$,  and  $B^{\ell}_{G}(v_{\ell +1})= \deg^+_{G}(v_{\ell +1})$ ,
let
$$B^0_{G-e}(v) =
    \begin{cases}
     B^0_{G}(v)-1 & \text{if }  v=v_a,\\
     B^0_{G}(v)+1 & \text{if }  v=v_{\ell+1},\\
        B^0_{G}(v) & \text{}  \mbox{otherwise.}

    \end{cases}$$
Then  for $t\in\{0,1,\ldots, a-1\}$,
$$B^t_{G-e}(v) =
    \begin{cases}
     B^t_{G}(v)-1 & \text{if }  v=v_a,\\
     B^0_{G}(v)+1 & \text{if }  v=v_{\ell+1},\\
        B^t_{G}(v) & \text{}  \mbox{otherwise.}
    \end{cases}$$
Now for $t\in\{a, a+1,\ldots,\ell\}$,
$$B^t_{G-e}(v) = B^t_{G}(v) \: \text{for all} \:v \in V(G) .$$
Observe that $B^{\ell}_{G-e}(v_{\ell+1})=B^{\ell}_{G}(v_{\ell+1})= \deg^+_{G-e}(v_{\ell+1})$  and so $v_{\ell+1}$ can fire. Consequently for $t\in\{\ell+1,\ell+2,\ldots, n\}$,
    $$B^t_{G-e}(v) = B^t_{G}(v) \: \text{for all} \:v \in V(G) .$$
    Thus  $B(G-e)\le B(G)$.

Case IV:   $v_a \in L \setminus \{v_{\ell}\}$ and  $v_b  \in L$.
\\
As the arc $e$ is not present  in $G-e$,  and  $B^{\ell}_{G}(v_{b})= \deg^+_{G}(v_{b})$
let
$$B^0_{G-e}(v) =
    \begin{cases}
     B^0_{G}(v)-1 & \text{if }  v=v_a,\\
     B^0_{G}(v)+1 & \text{if }  v=v_{b},\\
        B^0_{G}(v) & \text{}  \mbox{otherwise.}

    \end{cases}$$
Then  for $t\in\{0,1,\ldots, a-1\}$,
$$B^t_{G-e}(v) =
    \begin{cases}
     B^t_{G}(v)-1 & \text{if }  v=v_a,\\
     B^0_{G}(v)+1 & \text{if }  v=v_{b},\\
        B^t_{G}(v) & \text{}  \mbox{otherwise.}
    \end{cases}$$
Now for $t\in\{a, a+1,\ldots,b\}$
$$B^t_{G-e}(v) = B^t_{G}(v) \: \text{for all} \:v \in V(G) .$$
Observe that $B^{b-1}_{G-e}(v_{b})= \deg^+_{G-e}(v_{b})$  and so $v_{b}$ can fire. Consequently for $t\in\{b,b+1,\ldots, n\}$,
    $$B^t_{G-e}(v) = B^t_{G}(v) \: \text{for all} \:v \in V(G) .$$
Thus $B(G-e)\le B(G)$.

Case V:   $v_a \in R$ and  $ v_b \in R \setminus \{v_{\ell +1}\}$.
\\
Let $B^0_{G-e}(v)=B^0_{G}(v)$ for each $v\in V(G)$. Then  for $t\in\{0,1,\ldots, a-1\}$,  $B^t_{G-e}(v) = B^t_{G}(v)$. Now for $t\in\{a, a+1,\ldots, b-1\}$,
$$B^t_{G-e}(v) =
    \begin{cases}
     B^t_{G}(v)-1 & \text{if }  v=v_b,\\
        B^t_{G}(v) & \text{}  \mbox{otherwise.}
    \end{cases}$$
Observe that $B^{b-1}_{G-e}(v_b)=B^{b-1}_{G}(v_b)-1\ge \deg^+_{G-e}(v_b)$  and so $v_b$ can fire. Consequently for $t\in\{b,b+1,\ldots, n\}$,
$$B^t_{G-e}(v) =
    \begin{cases}
     B^t_{G}(v)-1 & \text{if }  v=v_b,\\
        B^t_{G}(v) & \text{}  \mbox{otherwise.}
    \end{cases}$$
    Thus  $B(G-e)\le B(G)$.

Case VI:  $v_a=v_{\ell}$, $v_b=v_{\ell+1}$ and $n$ is odd.
\\
As the arc $e$ is not present  in $G-e$,  and  $B^{\ell}_{G}(v_{\ell +1})= \deg^+_{G}(v_{\ell +1})$,
let
$$B^0_{G-e}(v) =
    \begin{cases}
     B^0_{G}(v)-1 & \text{if }  v=v_\ell,\\
     B^0_{G}(v)+1 & \text{if }  v=v_{\ell+1},\\
        B^0_{G}(v) & \text{}  \mbox{otherwise.}

    \end{cases}$$
Then  for $t\in\{0,1,\ldots, \ell-1\}$,
$$B^t_{G-e}(v) =
    \begin{cases}
     B^t_{G}(v)-1 & \text{if }  v=v_\ell,\\
     B^0_{G}(v)+1 & \text{if }  v=v_{\ell+1},\\
        B^t_{G}(v) & \text{}  \mbox{otherwise.}
    \end{cases}$$
Now for $t\in\{a, a+1,\ldots,\ell\}$,
$$B^t_{G-e}(v) = B^t_{G}(v) \: \text{for all} \:v \in V(G) .$$
Observe that $B^{\ell}_{G-e}(v_{\ell+1})= \deg^+_{G-e}(v_{\ell+1})$  and so $v_{\ell+1}$ can fire. Consequently for $t\in\{\ell+1,\ell+2,\ldots, n\}$,
    $$B^t_{G-e}(v) = B^t_{G}(v) \: \text{for all} \:v \in V(G) .$$
    Thus  $B(G-e)\le B(G)$.

Case VII:  $v_a=v_{\ell}$, $v_b=v_{\ell+1}$ and $n$ is even
\\
As the arc $e$ is not present  in $G-e$,
let $B^0_{G-e}(v_\ell)=B^0_{G}(v_\ell)-1$  and  $B^0_{G-e}(v)=B^0_{G}(v)$ for each $v\in V(G) \setminus \{v_\ell\}$.
Then  for $t\in\{0,1,\ldots, \ell-1\}$,
$$B^t_{G-e}(v) =
    \begin{cases}
     B^t_{G}(v)-1 & \text{if }  v=v_\ell,\\
        B^t_{G}(v) & \text{}  \mbox{otherwise.}
    \end{cases}$$
Now for $t=\ell$,
$$B^t_{G-e}(v) =
    \begin{cases}
     B^t_{G}(v)-1 & \text{if }  v=v_{\ell +1},\\
        B^t_{G}(v) & \text{}  \mbox{otherwise.}
    \end{cases}$$
Observe that $B^{\ell}_{G-e}(v_{\ell+1})=B^{\ell}_{G}(v_{\ell+1})-1= \deg^+_{G-e}(v_{\ell+1})$  and so $v_{\ell +1}$ can fire. Consequently for $t\in\{\ell+1,\ell+2,\ldots, n\}$,
$$B^t_{G-e}(v) =
    \begin{cases}
     B^t_{G}(v)-1 & \text{if }  v=v_{\ell +1},\\
        B^t_{G}(v) & \text{}  \mbox{otherwise.}
    \end{cases}$$
Thus $B(G-e)\le B(G)-1$.

Case VIII: $v_a=v_{\ell}$ and   $v_b\in R \setminus \{v_{\ell +1}\}$.
\\
As the arc $e$ is not present  in $G-e$,
let $B^0_{G-e}(v_\ell)=B^0_{G}(v_\ell)-1$  and  $B^0_{G-e}(v)=B^0_{G}(v)$ for each $v\in V(G) \setminus \{v_\ell\}$. Then  for $t\in\{0,1,\ldots, \ell-1\}$,
$$B^t_{G-e}(v) =
    \begin{cases}
     B^t_{G}(v)-1 & \text{if }  v=v_\ell,\\
        B^t_{G}(v) & \text{}  \mbox{otherwise.}
    \end{cases}$$
Now for $t\in\{\ell, \ell+1,\ldots, b-1\}$,
$$B^t_{G-e}(v) =
    \begin{cases}
     B^t_{G}(v)-1 & \text{if }  v=v_b,\\
        B^t_{G}(v) & \text{}  \mbox{otherwise.}
    \end{cases}$$
Observe that $B^{b-1}_{G-e}(v_b)=B^{b-1}_{G}(v_b)-1\ge \deg^+(v_b)$  and so $v_b$ can fire. Consequently for $t\in\{b,b+1,\ldots, n\}$,
$$B^t_{G-e}(v) =
    \begin{cases}
     B^t_{G}(v)-1 & \text{if }  v=v_b,\\
        B^t_{G}(v) & \text{}  \mbox{otherwise.}
    \end{cases}$$
    Thus  $B(G-e)\le B(G)-1$.
    \end{proof}
In contrast to this result for transitive tournaments, in general for directed acyclic graphs it is possible for the brushing number to increase when an arc is removed. For an example consider the graph $G$ in Figure~\ref{fig:10}. The brushing number of $G$ is 3. When  the arc $e$ is removed, the brushing number of $G-e$ is 5.
\usetikzlibrary{graphs}

\subsection{Directed Acyclic Graphs}\label{sec 2.2}
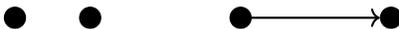
\begin{figure}[h]
    \centering
    \begin{tikzpicture}
          \node[normal] (v1) at (-8,0) {};
        \node[normal] (v2) at (-9,0) {};
          \node[normal] (v3) at (-6,0) {};
           \node[normal] (v4) at (-4,0) {};
         \draw[thick,->] (v3) to (v4);
    \end{tikzpicture}
    \caption{The set of directed acyclic graphs with two vertices.}
    \label{fig:1}
\end{figure}

Denote $\mathfrak{G}_k$ as the set of all directed acyclic graphs on $k$ vertices.  For each $k\ge2$  we will find upper bounds for the brushing number of all elements of $\mathfrak{G}_k$. For $k=2$, there are two directed acyclic graphs as shown in Figure \ref{fig:1} with brushing number 2 and 1 respectively. We will continue with the case where $k=3$.
\begin{lemma}\label{le_1.2}
    If $G\in\mathfrak{G}_3$ then $B(G)\leq 3 $.
\end{lemma}
\begin{proof}
    All graphs of $\mathfrak{G}_3$ are shown in Figure \ref{fig:G_4_5} and each can be brushed with at most three brushes.
    \end{proof}
    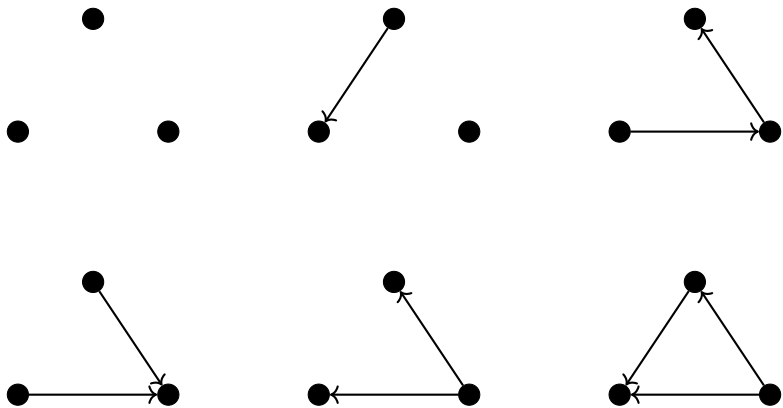
\begin{figure}[h]
    \centering
    \begin{tikzpicture}

          \node[normal] (v9) at (-8,0) {};
        \node[normal] (v10) at (-9,-1.5) {};
        \node[normal] (v11) at (-7,-1.5) {};

          \node[normal] (v13) at (-3,-1.5) {};
        \node[normal] (v14) at (-4,0) {};
        \node[normal] (v15) at (-5,-1.5) {};
         \draw[thick,->] (v14) to (v15);

  \node[normal] (v1) at (1,-1.5) {};
        \node[normal] (v2) at (0,0) {};
        \node[normal] (v3) at (-1,-1.5) {};
         \draw[thick,->] (v3) to (v1);
         \draw[thick,->] (v1) to (v2);

 \node[normal] (v4) at (-8,-3.5) {};
        \node[normal] (v5) at (-9,-5) {};
        \node[normal] (v6) at (-7,-5) {};
         \draw[thick,->] (v5) to (v6);
         \draw[thick,->] (v4) to (v6);

\node[normal] (v16) at (-3,-5) {};
        \node[normal] (v17) at (-4,-3.5) {};
        \node[normal] (v18) at (-5,-5) {};
         \draw[thick,->] (v16) to (v17);
         \draw[thick,->] (v16) to (v18);

\node[normal] (v19) at (1,-5) {};
        \node[normal] (v20) at (0,-3.5) {};
        \node[normal] (v21) at (-1,-5) {};
         \draw[thick,->] (v19) to (v20);
         \draw[thick,->] (v19) to (v21);
        \draw[thick,->] (v20) to (v21);

    \end{tikzpicture}
    \caption{All graphs of $\mathfrak{G}_3$.}
    \label{fig:G_4_5}
\end{figure}

The upper bound for the brushing number of graphs in $\mathfrak{G}_2$ and $\mathfrak{G}_3$ deviates from the general result obtained for the upper bound that we will establish for $k\ge 4$\ (Theorem~\ref{th:6}). Mathematical induction is going to be used in the proof of Theorem~\ref{th:6}, by considering the two cases where $k$ is even and odd. As the base cases of this proof we will give the upper bound for the brushing number of the two sets $\mathfrak{G}_4$ and $\mathfrak{G}_5$ as separate lemmata before proving Theorem~\ref{th:6}. In the following proofs, $S_{n,m}\subseteq\mathfrak{G}_n$ will denote the set of directed acyclic graphs with $n$ vertices and $m$ components.

\begin{lemma} \label{lem:G_3}
    If $G\in\mathfrak{G}_4$ then $B(G)\leq 4$.
\end{lemma}

\begin{proof}
There is one graph $G_1\in S_{4,4}\subset \mathfrak{G}_4$ with four components which is shown in Figure~\ref{fig:G_4_1}. The graph $G_1$ can brushed with four brushes.
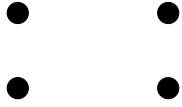
\begin{figure}[h]
    \centering
    \begin{tikzpicture}

        \node[normal] (v1) at (0,0) {};
        \node[normal] (v2) at (0,-1) {};
        \node[normal] (v3) at (2,0) {};
         \node[normal] (v13) at (2,-1) {};
    \end{tikzpicture}
    \caption{$G_1\in \mathfrak{G}_4$ with four components.}
    \label{fig:G_4_1}
\end{figure}
Also there is one graph $G_2\in S_{4,3}\subset\mathfrak{G}_4$ with three components which is shown in Figure~\ref{fig:G_4_2}. The graph $G_2$ can brushed with three brushes.
\begin{figure}[h]
    \centering
    \begin{tikzpicture}

        \node[normal] (v1) at (0,0) {};
        \node[normal] (v2) at (0,-1) {};
        \node[normal] (v3) at (2,0) {};
         \node[normal] (v13) at (2,-1) {};
          \draw[thick,->] (v1) to (v3);
    \end{tikzpicture}
    \caption{$G_2 \in \mathfrak{G}_4$ with three components.}
    \label{fig:G_4_2}
\end{figure}
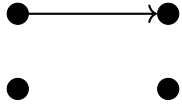
The graphs of $S_{4,2}\subset \mathfrak{G_4}$ are shown in Figure~\ref{fig:G_4_3} and each can be brushed with at most three brushes.

\begin{figure}[h]

    \centering
    \begin{tikzpicture}

        \node[normal] (v1) at (0,0) {};
        \node[normal] (v2) at (0,-1) {};
        \node[normal] (v3) at (2,0) {};
         \node[normal] (v4) at (2,-1) {};
        \draw[thick,->] (v1) to (v3);
         \draw[thick,->] (v2) to (v4);

        \node[normal] (v5) at (-4,0) {};
        \node[normal] (v6) at (-4,-1) {};
        \node[normal] (v7) at (-2,0) {};
         \node[normal] (v8) at (-2,-1) {};
        \draw[thick,->] (v6) to (v8);
         \draw[thick,->] (v8) to (v7);

          \node[normal] (v9) at (-8,0) {};
        \node[normal] (v10) at (-8,-1) {};
        \node[normal] (v11) at (-6,0) {};
         \node[normal] (v12) at (-6,-1) {};
        \draw[thick,->] (v10) to (v12);
         \draw[thick,->] (v11) to (v12);

          \node[normal] (v13) at (-8,-3) {};
        \node[normal] (v14) at (-8,-4) {};
        \node[normal] (v15) at (-6,-3) {};
         \node[normal] (v16) at (-6,-4) {};
        \draw[thick,->] (v16) to (v14);
         \draw[thick,->] (v16) to (v15);

        \node[normal] (v17) at (-4,-3) {};
        \node[normal] (v18) at (-4,-4) {};
        \node[normal] (v19) at (-2,-3) {};
         \node[normal] (v20) at (-2,-4) {};
        \draw[thick,->] (v20) to (v18);
         \draw[thick,->] (v20) to (v19);
           \draw[thick,->] (v18) to (v19);

\begin{comment}

        \node[normal] (v4) at (0,1) {};
        \node[normal] (v5) at (1,1) {};
        \node[normal] (v6) at (2,1) {};
        \draw[thick,->] (v4) to (v5);
        \draw[thick,->] (v4) to [out=45,in=135](v6);

        \node[normal] (v7) at (3,1) {};
        \node[normal] (v8) at (4,1) {};
        \node[normal] (v9) at (5,1) {};
        \draw[thick,->] (v7) to (v8);
        \draw[thick,->] (v8) to (v9);
        \draw[thick,->] (v7) to [out=45,in=135](v9);

        \node[normal] (v10) at (3,0) {};
        \node[normal] (v11) at (4,0) {};
        \node[normal] (v12) at (5,0) {};
        \draw[thick,->] (v11) to (v12);
        \draw[thick,->] (v10) to [out=45,in=135](v12);
\end{comment}
  \end{tikzpicture}
    \caption{The graphs $S_{4,2}\subset \mathfrak{G_4}$.}
    \label{fig:G_4_3}
\end{figure}
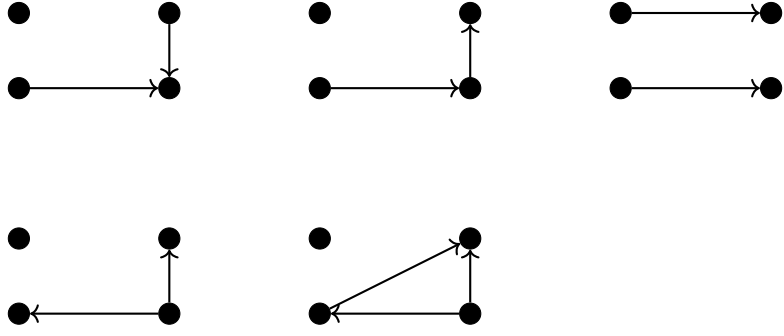

Now we will consider the set of graphs $S_{4,1}\subset \mathfrak{G}_4$. For the graphs in $S_{4,1}$ we will calculate the brushing number of graphs with $6, 5, 3$ and $4$ edges respectively. In $S_{4,1}$ the graphs with six edges are transitive tournaments. The brushing number of a transitive tournament with four vertices is $4$ as proved in Theorem~\ref{th_1.3}. Consequently from Theorem~\ref{th_1.4}, it follows that when $G$ is a directed acyclic graph with four vertices, five edges and one component $B(G)\leq 4$.
The set of graphs from $S_{4,1}$ with three edges is   illustrated in Figure~\ref{fig:G_4_4}.

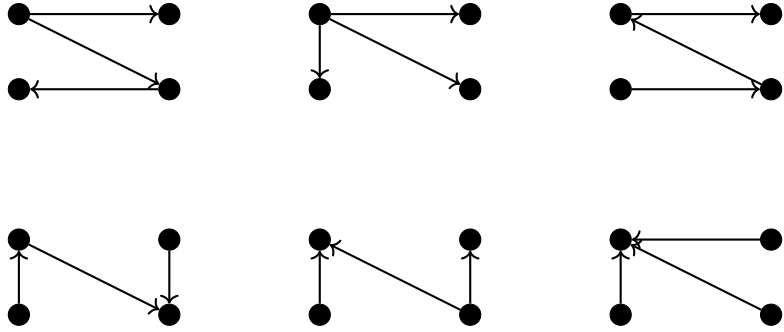
\begin{figure}[h]
    \centering
    \begin{tikzpicture}

        \node[normal] (v1) at (0,0) {};
        \node[normal] (v2) at (0,-1) {};
        \node[normal] (v3) at (2,0) {};
         \node[normal] (v4) at (2,-1) {};
        \draw[thick,->] (v1) to (v3);
         \draw[thick,->] (v2) to (v4);
          \draw[thick,->] (v4) to (v1);

        \node[normal] (v5) at (-4,0) {};
        \node[normal] (v6) at (-4,-1) {};
        \node[normal] (v7) at (-2,0) {};
         \node[normal] (v8) at (-2,-1) {};
        \draw[thick,->] (v5) to (v7);
         \draw[thick,->] (v5) to (v8);
          \draw[thick,->] (v5) to (v6);

          \node[normal] (v9) at (-8,0) {};
        \node[normal] (v10) at (-8,-1) {};
        \node[normal] (v11) at (-6,0) {};
         \node[normal] (v12) at (-6,-1) {};
        \draw[thick,->] (v9) to (v11);
         \draw[thick,->] (v9) to (v12);
         \draw[thick,->] (v12) to (v10);

          \node[normal] (v13) at (-8,-3) {};
        \node[normal] (v14) at (-8,-4) {};
        \node[normal] (v15) at (-6,-3) {};
         \node[normal] (v16) at (-6,-4) {};
        \draw[thick,->] (v14) to (v13);
         \draw[thick,->] (v15) to (v16);
         \draw[thick,->] (v13) to (v16);

        \node[normal] (v17) at (-4,-3) {};
        \node[normal] (v18) at (-4,-4) {};
        \node[normal] (v19) at (-2,-3) {};
         \node[normal] (v20) at (-2,-4) {};
        \draw[thick,->] (v18) to (v17);
         \draw[thick,->] (v20) to (v17);
           \draw[thick,->] (v20) to (v19);

             \node[normal] (v21) at (0,-3) {};
        \node[normal] (v22) at (0,-4) {};
        \node[normal] (v23) at (2,-3) {};
         \node[normal] (v24) at (2,-4) {};
        \draw[thick,->] (v22) to (v21);
         \draw[thick,->] (v23) to (v21);
          \draw[thick,->] (v24) to (v21);

    \end{tikzpicture}
    \caption{The set of graphs from $S_{4,1}$ with three edges.}
    \label{fig:G_4_4}
\end{figure}

All graphs in Figure \ref{fig:G_4_4} can be brushed with  three brushes. Now consider the set of graphs from $S_{4,1}$ with four edges. Let $G \in S_{4,1}$ with four edges. Each graph $G$ has a spanning sub-tree $H$ with three edges. If $e\in A(G)$ and $e \notin A(H)$, then $ G-e$ is connected. As $G-e = H$ will be one of the graphs in $S_{4,1}$ with three edges, which can be brushed using at most three brushes, the brushing number of $G \in S_4$ is at most $4$.
\end{proof}

\begin{lemma}\label{lemma_2.3}
    If $G\in\mathfrak{G}_5$ then $B(G)\leq 6$.
\end{lemma}

\begin{proof}
Consider a graph $G\in \mathfrak{G}_{5}$. Pick a source $u$ and a sink $v$ from $G$ and consider the remaining set of vertices, $S=V(G)\setminus\{u,v\}$. Note that $|S|=3$. Let $X$, $Y$ and $Z$ be subsets of $S$ such that
$X=N^+(u)\cap N^-(v)$ (i.e., $X$ consists of vertices that are each an out-neighbour of $u$ and an in-neighbour of $v$),  $Y=\{N^+(u)\cup N^-(v)\}\setminus X$, and $Z=S\setminus\{X \cup Y\}$. Note $|X|+|Y|+|Z|= 3$, so
%Also $0\le |X|\le 3$,  $0\le |Y|\le 3$, and $0\le |Z|\le 3$.
the graph $G-\{u,v\}$ is a directed acyclic graph on three vertices.

First consider the case where $G-\{u,v\}$ has less than three components. In $G$ place $|X|$ brushes at $u$, to later brush paths of the form $(u,x,v)$ where $x\in X$.
Partition $Y$ so that $Y_1 \subseteq Y$ is the set of vertices with incident arcs of the form $(u,y)$ where $y\in Y_1$ and $Y_2 \subseteq Y$ is the set of vertices with incident arcs of the form $(y,v)$ where $y\in Y_2$. Note that $|Y|=|Y_1|+|Y_2|$.
%Let $Y_1 \subseteq Y$ be the set of vertices with incident arcs of the form $(u,y)$ where $y\in Y_1$. Let $Y_2 \subseteq Y$ be the set of vertices with incident arcs of the form $(y,v)$ where $y\in Y_2$. Also $|Y|=|Y_1|+|Y_2|$.
Then place $|Y_1|$ brushes at $u$ to brush arcs of the form $(u,y)$ and place $|Y_2|$ brushes, one brush at each vertex of $Y_2$, to brush arcs of the form $(y,v)$ where $y\in Y_2$. Also, if there is an arc from $u$ to $v$, place one more brush at $u$.
However, if $u$ is isolated in $G$ then place one brush at $u$.
Similarly, if $v$ is isolated in $G$ then place a brush at $v$.
As proved in Lemma~\ref{le_1.2} we know that when $G-\{u,v\}$ has less than three components the brushing number of $G-\{u,v\}$ is at most $2$. Then place up to $2$ brushes on set $S$ in a way such that it will be consistent with a  strategy that successfully brushes $G-\{u,v\}$. Now first fire vertex $u$ and then all the vertices in $S$ in accordance with the strategy used to brush $G-\{u,v\}$. After all vertices in $S$ fire, $v$ accumulates brushes as $v$ is the sink.  Finally vertex $v$ fires.
If both $u$ and $v$ are isolated in $G$, then $B(G) = 2 + B(G-\{u,v\}) \le 4$; otherwise
%If $u$ and $v$ are not both isolated in $G$, then
$$B(G)\le |X|+|Y|+1+2\le 3+1+2 =6\mbox{.}$$

When $G-\{u,v\}$ has three components, place $|X|+|Y_1|$ brushes at $u$ and also one brush at each vertex of $Y_2$. As $G-\{u,v\}$ has no arcs now, place $|Z|$ brushes, one brush each at every vertex in $Z$. Then fire the vertices in a sequence similar to the previous case. Therefore, $$B(G)\le |X|+|Y|+|Z|+2\le 3+2=5.$$
\end{proof}

Having established upper bounds for the brushing number of $G\in\mathfrak{G}_k$, for small values of $k$, we now consider larger values of $k$.

\begin{theorem}\label{th:6}
If $k\ge4$ and $G\in\mathfrak{G}_k$, then

$$B(G)\le
    \begin{cases}
     \frac{k^2-1}{4} & \text{if } \mbox{k is odd,}\\
        \frac{k^2}{4} & \text{if } \mbox{k is even.}
    \end{cases}$$
\end{theorem}

\begin{proof}
%The method of mathematical induction is used in this proof.
Note that the theorem holds when $k \in \{4,5\}$, as established by Lemma~\ref{lem:G_3} and Lemma~\ref{lemma_2.3}.
To complete the proof, we now show that if the theorem holds for every directed graph in $\mathfrak{G}_n$ then it also holds for every directed graph in $\mathfrak{G}_{n+2}$.

Let $G\in \mathfrak{G}_{n+2}$. Pick a source $u$ and a sink $v$ from $G$ and consider the remaining set of vertices, $S=V(G)\setminus\{u,v\}$. Note that $|S|=n$ and so the subgraph $G'=G-\{u,v\}$ is in $\mathfrak{G}_n$. Let $X$, $Y$, and $Z$ be subsets of $S$ such that $X=N^+(u)\cap N^-(v)$,  $Y=\{N^+(u)\cup N^-(v)\}\setminus X$, and $Z=S\setminus\{X \cup Y\}$; observe that $|X|+|Y|+|Z|= n$.
%Also $0\le |X|\le n$,  $0\le |Y|\le n$, and $0\le |Z|\le n$.

Since the graph $G'$ is in $\mathfrak{G}_n$, then $B(G') \leq \lfloor \frac{n^2}{4} \rfloor$.  We will use a brushing strategy for $G'$ to devise a strategy for $G$.  Accordingly, distribute $B(G')$ brushes among the vertices of $S$ in $G$ so that they adhere to the initial configuration of a brushing strategy for $G'$.  We now place additional brushes in $G$ as follows.

Place $|X|$ brushes at $u$, to later brush paths of the form $(u,x,v)$ where $x\in X$. Partition $Y$ so that $Y_1 \subseteq Y$ is the set of vertices with incident arcs of the form $(u,y)$ where $y\in Y_1$ and $Y_2 \subseteq Y$ is the set of vertices with incident arcs of the form $(y,v)$ where $y\in Y_2$. Note that $|Y|=|Y_1|+|Y_2|$.  Place $|Y_1|$ additional brushes at $u$, to brush arcs of the form $(u,y)$, and one brush at each vertex of $Y_2$, to brush arcs of the form $(y,v)$ where $y\in Y_2$.
If there is an arc from $u$ to $v$ in $G$ then place one more brush at $u$.
However, if $u$ is isolated in $G$ then place one brush at $u$.
Similarly, if $v$ is isolated in $G$ then place a brush at $v$.

It follows that a successful brushing strategy is to fire $u$, then fire the vertices of $S$ in the same sequence as the strategy used for cleaning $G'$, and finally fire $v$.
Unless both $u$ and $v$ are isolated vertices in $G$, then
$B(G)\le B(G') + |X|+|Y|+1 \le \lfloor\frac{n^2}{4}\rfloor + n + 1$.
%and the theorem's conclusion easily follows.
In the case that $u$ and $v$ are both isolated, then we have $B(G) = B(G') + 2 < \lfloor\frac{n^2}{4}\rfloor + 2$.
In either case, the theorem's conclusion is easily obtained.
\end{proof}

\section{Other Types of Directed Graphs}\label{chap 3}

In this section we will find the brushing number of several types of directed graphs such as complete directed graphs, rooted trees, and rotational tournaments. Note that a \textit{complete directed graph} is a simple directed graph such that between each pair of its vertices, both (oppositely directed) arcs exist.
\begin{theorem}
If $G$ is a complete directed graph, then $$B(G)=\frac{|A(G)|}{2}.$$
\end{theorem}

\begin{proof}
Every vertex of a complete directed graph $G$ with $n$ vertices has out-degree $n-1$. Therefore, the total number of arcs $|A(G)|= n(n-1)$.
%(because every arc of $G$ has one tail and one head.)
Consider a complete directed graph $G$ with $n$ vertices, $v_1,v_2,\ldots,v_n$. The vertices of $G$ are labelled arbitrarily $v_1,v_2,\ldots,v_n$. Recall that $B^t_{G}(v)$ is the number of brushes at vertex $v$ of $G$, at time $t$ ($t=0,1,2, \ldots, n$). The vertex $v_i$ fires in between time $i-1$ and $i$. Consider the initial configuration where
$$B^0_{G}(v_{i}) = n-i \;\; \text{for }  i=1,2,\ldots, n.$$
As vertex $v_1$ with out-degree $n-1$ has  $n-1$  brushes,  $v_1$ fires.  Now
$$B^1_{G}(v_{i}) =
    \begin{cases}
    0 & \text{if } i=1, \\
     n-i+1 & \text{otherwise.}
    \end{cases}$$
Observe that the number of brushes at $v_2$ is $n-1$. So $v_2$ fires.  Then
$$B^2_{G}(v_{i}) =
    \begin{cases}
    1 & \text{if }  i=1,\\
     0
      & \text{if }  i=2,\\
        n-i+2 & \text{otherwise.}
    \end{cases}$$
After the $k^{\text{th}}$ vertex fires,
$$B^k_{G}(v_{i}) =
    \begin{cases}
    k-i
      & \text{if }  i=1, 2,\ldots, k-1,\\
     0
      & \text{if }  i=k\\
        n-i+k & \text{}  \mbox{otherwise.}
    \end{cases}$$
When $\{v_1, v_2,\ldots,v_{k-1}\}$ have all fired $v_k$ is able to fire, as it will have $n-k+k-1=n-1=\deg^+(v_k)$ brushes, where $2\le k\le n$. Similarly, we can also show that $v_{k+1},v_{k+2},\ldots, v_n$ are able to fire in sequence. The total number of brushes used is equal to
$(n-1)+ (n-2)+(n-3)+\cdots+2+1+0 = n(n-1)/2 = \frac{|A(G)|}{2}$ and hence $B(G)\le\frac{|A(G)|}{2}$.

For each $i=1,2,\ldots,n$ let $m_i$ be the number of in-neighbour vertices of $v_i$ that have fired before the vertex $v_i$ fires. For a vertex $v_i$ to fire, it is necessary that
%For each $i=1,2,\ldots,n$ let $m_i$ be the number of in-neighbour vertices fired before firing the $i^{\text{th}}$ vertex. For vertex $v_i$ to fire
\begin{align}
\label{inq2}
    B^0_{G}(v_{i})+m_i\ge (n-1).
\end{align}
Summing the inequality (\ref{inq2}) for all $i\in \{1,2, \ldots, n\}$ and observing that $m_i=i-1$ for $i\in \{1,2, \ldots, n\}$ yields
$$\sum_{i=1}^{n} B^0_{G}(v_{i})+\sum_{i=1}^{n}m_i\ge  n(n-1),$$
where $B^0_{G}(v_{i})$ is any valid initial brushing configuration that successfully cleans the graph. If $B^0_{G}(v_{i})$ is an optimal brushing configuration then $\sum_{i=1}^{n} B^0_{G}(v_{i})=B(G).$  Hence
$$B(G)\ge  {n(n-1)}-\sum_{i=1}^{n}(i-1)$$ and
$$B(G)\ge n(n-1)-\frac{(n^2-3n)}{2}.$$ Therefore $$B(G) \ge \frac{n(n-1)}{2}=\frac{|A(G)|}{2}.$$
\end{proof}
A \emph{directed tree} is defined as a directed graph whose underlying undirected graph is a tree. A \emph{rooted tree} is a directed tree with a distinguished vertex $r$, called the \emph{root} such that for every other vertex $v$, the unique path from $r$ to $v$ is a directed path from $r$ to $v$.
\begin{theorem}\label{th 6}
    If $G$ is a rooted tree with $k$ leaves, then $B(G)=k$.
\end{theorem}
\begin{proof}
Let $L\subseteq V(G)$ be the set of leaves of $G$, where $|L|=k$. Since there exist $k$ sinks,  $B(G)\ge k$ (by Lemma~\ref{le_1.5}).

There is a unique directed path from the root $r$ to each leaf $\ell\in L$. The union of the arcs of these paths is the arc set of $G$. Now place $k$ brushes at the root $r$ of the graph $G$. To clean the graph $G$ one brush can travel in each path from $r$ to each $\ell\in L$. Therefore, $B(G)\le k$.
\end{proof}
A \emph{path decomposition} of a directed graph $G$ is a set of arc-disjoint directed paths such that the union of the arcs of these paths is the arc set of $G$. If $M$ is a path decomposition of a directed graph $G$, and there does not exist a path decomposition of $G$ with less than $|M|$ paths, then $M$ is a minimum path decomposition of $G$.

\begin{theorem}\label{th 3}
    If $M$ is a minimum path decomposition of a directed acyclic graph $G$ then $B(G)\le |M|$.
\end{theorem}

\begin{proof}
Obtain a minimum path decomposition $M$ of a directed acyclic graph $G$. The paths of $G$ that we consider in this proof are elements of $M$. Let $a_i$ be the number of paths of $M$ that begin at  $v_i\in V(G)$ where $V(G)=\{v_1,v_2,\ldots,v_{|V(G)|}\}$. For all $v_i\in V(G)$ let $B^0_{G}(v_i)=a_i$. Let $\omega$ be the length of the longest path in $G$. For each $k\in \{0,1, \ldots,\omega\}$, let $S_k$ be the set of vertices with distance $k$ from the furthest source vertex in $G$  and let $S_{<k}=S_0\cup S_1\cup\cdots \cup S_{k-1}$. As $S_0$ is the set of source vertices in $G$ then for $v_i\in S_0, \; \deg^+(v_i)=a_i.$
We depart from our usual convention of sequentially firing vertices
and allow all vertices in $S_0$ to fire in parallel between time steps $0$ and $1$; the resulting configuration corresponds to that which would be obtained by sequentially firing the vertices in arbitrary order.

Consider $v_i\in S_k$ where $k\in \{1,2, \ldots,\omega\}$. Note that vertex $v_i$ receives brushes from its  in-neighbours in $S_{<k}$ before the $k^{\text{th}}$ time step. Observe that $\deg^+(v_i)$ is equal to the sum of the number of paths in $M$ that start at $v_i$ and the number of paths in $M$ that include $v_i$ as a middle vertex. Hence $\deg^+(v_i)$ is at most the sum of  $a_i$ and the number of brushes that $v_i$ receives from its in-neighbours in $S_{<k}$, which is equal to
$B^k_{G}(v_i)$. Therefore $\deg^+(v_i)\le B^k_{G}(v_i)$. Then each vertex of $S_k$ fires (in parallel) between the time steps $k$ and $k+1$ where $k\in \{1,2, \ldots,\omega\}$. Hence $B(G)\le\sum\limits_{v\in V(G)}B^0_{G}(v)=| M|$.
\end{proof}
The \emph{transpose} of a directed graph $G$ is another directed graph $G^T$ on the same set of vertices with all the arcs reversed compared to the orientation of the corresponding arcs in $G$.
\begin{lemma}\label{le 1}
    If $G$ is a directed acyclic graph, then $B(G^T)\le B(G).$
\end{lemma}
\begin{proof}
Let $G$ be a directed acyclic graph with the vertex set $\{v_1,v_2,\ldots,v_n\}$ and optimal initial brushing configuration $B^0_{G}(v_i)=a_i$ for all $ v_i\in V(G)$ that successfully cleans the graph $G$. Observe that at some time $t'$ the graph $G$ is clean and the brushing configuration is $B^{t'}_{G}(v_i)=b_i,\;\text{for all}\;v_i\in V(G)$. Throughout the cleaning process of $G$ each brush in the initial brushing configuration $B^0_{G}(v_i)$ traverses  a directed path. If the transpose of these directed paths is considered then we obtain a path decomposition $M'$ of  $G^T$. Define the initial brushing configuration of $G^T$ as $B^0_{G^T}(v_i)=B^{t'}_{G}(v_i)=b_i$.  Then  $b_i$ is the number of paths of $M'$ that begin at  $v_i$ in $G^T$.

Now we use a similar approach as in Theorem~\ref{th 3}. Let $\omega$ be the length of the longest path in $M'$. For each $k\in \{0,1, \ldots,\omega\}$ let $S_k$ be the set of vertices with distance $k$ from the furthest source vertex in $G^T$  and $S_{<k}=S_0\cup S_1\cup\cdots \cup S_{k-1}$. As $S_0$ is the set of source vertices in  $G^T$ for $v_i\in S_0, \; \deg^+_{G^T}(v_i)=b_i.$   Fire each vertex of $S_0$ (in parallel) between time steps $0$ and $1$.

Consider $v_i\in S_k$ where $k\in \{1,2, \ldots,\omega\}$. Note that vertex $v_i$ receives brushes from its  in-neighbours in $S_{<k}$ before the $k^{\text{th}}$ time step. Observe that $\deg^+_{G^T}(v_i)$ is equal to the sum of the number of paths in $M'$ that start at $v_i$ and the number of paths in $M'$ that include $v_i$ as a middle vertex. Hence $\deg^+_{G^T}(v_i)$ is at most the sum of  $b_i$ and the number of brushes that $v_i$ receives from its in-neighbours in $S_{<k}$, which is equal to
$B^k_{G}(v_i)$. Therefore $\deg^+_{G^T}(v_i)\le B^k_{G}(v_i)$. Then each vertex of $S_k$ fires (in parallel) between the time steps $k$ and $k+1$ where $k\in \{1,2, \ldots,\omega\}$. Hence $B(G^T)\le\sum\limits_{v\in V(G^T)}B^0_{G^T}(v)=\sum\limits_{v\in V(G)}B^{t_1}_{G}(v)=\sum\limits_{v\in V(G)}B^0_{G}(v)=B(G)$.
\end{proof}

\begin{theorem}
    If $G$ is a directed acyclic graph, then $B(G^T)=B(G)$.
\end{theorem}

\begin{proof}
    From from Lemma~\ref{le 1} we have $B(G^T)\le B(G)$. We also have $B(G)=B((G^T)^T)\le B(G^T)$.
\end{proof}
The following principle is described in \cite{bondy2008graph}.
\newline
\textbf{Principle of Directional Duality}: Any statement about a directed graph has an accompanying dual statement, obtained by applying the statement to the transpose of the directed graph and reinterpreting it in terms of the original directed graph.
\begin{theorem}\label{th10}
    If $T$ is a directed tree on $n$ vertices, the set of source vertices of $T$ is $S_0$ and the set of sink vertices is $S_\infty$, then $$B(T)\ge \max\left\{ \sum_{v\in S_0} |N^+(v)|, \sum_{v\in S_\infty} |N^-(v)|\right \}.$$
\end{theorem}
\begin{proof}
    Each arc incident with a vertex of $S_0$ requires a distinct brush (no brush can clean two of the incident arcs). Therefore $B(T)\ge \sum\limits_{v\in S_0} |N^+(v)|$. By the principle of directional duality $B(T)\ge \sum\limits_{v\in S_\infty} |N^-(v)|$.
\end{proof}

A \textit{rotational tournament} is defined in \cite{gross2003handbook} as follows. Let $\Gamma$ be an abelian group of odd order $n = 2m+1$ with identity $0$. Let $S$ be an $m$-element subset of $\Gamma \setminus \{0\}$ such that for every $x, y \in S,\, x + y \neq 0$. That is, choose exactly one element from each of the $m$ $2$-sets of the form $\{x, −x\}$, where $x$ ranges over all $x \in \Gamma \setminus \{0\}$. Form the directed graph $D$ with vertex set $V(D) = \Gamma$ and arc set $A(D)$ defined by: arc $(x, y) \in A(D)$ if and only if $y − x \in S$. Then $D$ is called a rotational tournament with symbol set $S$ and is denoted $R_\Gamma(S)$, or simply $R(S)$ if the group $\Gamma$ is understood.

A \textit{regular tournament} is a tournament $T$ in which there is an integer $s$ so that $\deg^+(v)=s$ for all vertices $v\in T$. The rotational tournament $R_\Gamma(S)$, where $|\Gamma|=n$, is a regular tournament on $n$ vertices. Each vertex in a rotational tournament with $n$ vertices has out-degree $\frac{n-1}{2}$.

\begin{theorem}\label{th 2}
    If $G$ is a rotational tournament with $n$ vertices then $B(G) = \frac{n^2-1}{8}.$
\end{theorem}

\begin{proof}
Let $B^t_{G}(v)$ be the number of brushes at vertex $v$ of $G$, at time $t$ ($t=0,1,2, \ldots, n$). The vertex $v_i$ fires in between time $i-1$ and $i$.  Let $B^0_{G}(v_{k})$ denote the number of brushes at the $k^{\text{th}}$ vertex in the initial configuration.

$$B^0_{G}(v_{k}) =
    \begin{cases}
     \frac{n-1}{2}-(k-1) & \text{if }  k=1, 2,\ldots,\frac{n-1}{2}\\
        0 & \text{}  \mbox{otherwise.}
    \end{cases}$$

So no vertex with index greater than $ \frac{n-1}{2}$ receives a brush in the initial configuration.  Recall that $\deg^+(v_i)= \frac{n-1}{2},\:\mbox{for all}\;0\le i\le n$. As vertex $v_1$ with out-degree $ \frac{n-1}{2}$ has  $ \frac{n-1}{2}$  brushes,  $v_1$ fires.  Now
$$B^1_{G}(v_{k}) =
    \begin{cases}

     \frac{n-1}{2}-(k-1)+1 & \text{if }  k=2,3,\ldots
     ,\frac{n-1}{2}\\
      1 & \text{if } k=\frac{n-1}{2} +1 \\
        0 & \text{}  \mbox{otherwise.}
    \end{cases}$$
Observe that the number of brushes at $v_2$ is $\frac{n-1}{2}$. So $v_2$ fires. Then

$$B^2_{G}(v_{k}) =
    \begin{cases}

     \frac{n-1}{2}-(k-1)+2 & \text{if }  k=3,4,\ldots,\frac{n-1}{2}\\
      2 & \text{if } k=\frac{n-1}{2} +1 \\
      1 & \text{if } k=\frac{n-1}{2} +2 \\
        0 & \text{}  \mbox{otherwise.}
    \end{cases}$$
Thus after $v_j$ vertex fires, where $1\le j\le \frac{n-1}{2}$.
$$B^j_{G}(v_{k}) =
    \begin{cases}

     \frac{n-1}{2}-(k-1)+j & \text{if }  k=j+1,j+2,\ldots,\frac{n-1}{2}\\
      j-(\ell-1) & \text{if } k=\frac{n-1}{2} +\ell ,(\ell=1,2,\ldots,j)\\
      0 & \text{}  \mbox{otherwise.}
    \end{cases}$$

When $\{v_1, v_2,\ldots,v_{j}\}$ have all fired $v_{j+1}$ is able to fire, as it will have $\frac{n-1}{2}-(j+1-1)+j-1=\frac{n-1}{2}=\deg^+(v_{j+1})$ brushes, where $1< j+1\le \frac{n-1}{2}$. After $v_{\frac{n-1}{2}+1}$ fires,
$$B^{\frac{n-1}{2}+1}_{G}(v_{k}) =
    \begin{cases}
    \frac{n-1}{2}-(\ell-2) & \text{if } k=\frac{n-1}{2} +\ell ,(\ell=2,3,\ldots,\frac{n+1}{2})\\
        0 & \text{}  \mbox{otherwise.}
    \end{cases}$$

After $v_{\frac{n-1}{2}+r}$ fires, where $2\le r\le \frac{n+1}{2}$,

$$B^{\frac{n-1}{2}+r}_{G}(v_{k}) =
    \begin{cases}
    \frac{n-1}{2}-[\ell-(r+1)] & \text{if } k=\frac{n-1}{2} +\ell ,(\ell=r+1,r+2,\ldots,\frac{n+1}{2})\\
    r-k & \text{if } k=1,2,\ldots,r-1\\
        0 & \text{}  \mbox{otherwise.}
    \end{cases}$$
Similarly, we can  also show that all $v_{\frac{n-1}{2}+r}$, where $2\le r\le \frac{n+1}{2}$, are able to fire in sequence.

When the above brushing strategy is used for $G$ we obtain,
$$ B(G)\le\sum_{k=1}^{ \frac{n-1}{2}} B^0_{G}(v_{k})=\frac{n-1}{2}+\frac{n-1}{2}-1+\frac{n-1}{2}-2+\cdots+1=\frac{n^2-1}{8}.$$

Let $m_i$ be the number of in-neighbour vertices of $v_i$ that have fired before the vertex $v_i$ fires. For a vertex $v_k$ to fire, it is necessary that
\begin{align}
\label{inq1}
    B^0_{G}(v_{k})+m_k\ge \frac{n-1}{2}
\end{align}
and
$$m_k =
    \begin{cases}
     k-1 & \text{if }  k=1, 2,\ldots, \frac{n-1}{2}\\
        \frac{n-1}{2} & \text{}  \mbox{otherwise.}
    \end{cases}$$
Summing the inequality (\ref{inq1}) for all $k\in \{1,2, \ldots, n\}$ we get
$$\sum_{k=1}^{n} B^0_{G}(v_{k})+\sum_{k=1}^{n}m_k\ge  \frac{n(n-1)}{2},$$
$$B(G)\ge  \frac{n(n-1)}{2}-\left(\sum_{k=1}^{\frac{n(n-1)}{2}}(k-1)+\sum_{k=\frac{(n+1)}{2}}^{n}\frac{(n-1)}{2}\right )$$ and
$$B(G)\ge  \frac{n(n-1)}{2}-\left(\frac{(n-3)}{2}\frac{(n-1)}{2}\frac{1}{2}+\frac{(n-1)}{2}\frac{(n+1)}{2}\right).$$ Therefore $$B(G) \ge \frac{n^2-1}{8}.$$

\end{proof}

\begin{corollary}
If $G$ is a rotational tournament that is cleaned using the strategy given in Theorem~\ref{th 2} then there is at least one brush that travels a Hamiltonian path during brushing $G$.

\end{corollary}
\begin{proof}
    When using the strategy given in Theorem~\ref{th 2} to brush a rotational tournament $G$, every vertex obtains  at least one brush from the preceding vertex that just fired.
\end{proof}

If $M$ is a minimum path decomposition of a directed graph $G$, then $| M |$ is called the \emph{path number} of $G$ and following the notation of \cite{https://doi.org/10.1112/plms.12328} it is denoted by $pn(G)$. A lower bound on $pn(G)$ can be given by considering the degree sequence of $G$. For each $v\in V(G)$ let $d_{v}=\deg^+_{G}(v)-\deg^-_{G}(v)$. Observe that in any path decomposition of $G$ at least $d_{v}$ paths should start at $v$. Hence $pn(G)\ge\frac{1}{2}\sum\limits_{v\in V(G)}|d_{v}|$.

In \cite{10.1007/978-3-030-83823-2_112} a \emph{perfect decomposition} of a directed graph $G$ is introduced as a set $P = \{P_1,\ldots,P_r\}$ of arc disjoint paths of $G$ that together cover $E(G)$, where $r = \frac{1}{2}\sum\limits_{v\in V(G)}|d_{v}|$.

\begin{theorem}
    \cite{https://doi.org/10.1112/plms.12328}
    If $G$ is a directed acyclic graph, then it has a perfect decomposition.
\end{theorem}

\begin{theorem}\label{th 4}
    If $G$ is a directed acyclic graph with a perfect decomposition $P = \{P_1,\ldots,P_r\}$, then $B(G)\le r$.
\end{theorem}

\begin{proof}
The strategy used in Theorem~\ref{th 3} can be generalised for any path decomposition. Therefore by following an approach similar to that of Theorem~\ref{th 3} we obtain $B(G)\le r$.
\end{proof}
Comparing Theorem~\ref{th 3} and Theorem~\ref{th 4} we obtain a better upper bound for the brushing number $B(G)$ of a directed acyclic graph which is: $pn(G)\ge r\ge B(G)$.

\section{Discussion}\label{sec 5}

\iffalse
It remains as an open problem to find an example of a directed acyclic graph with the property $pn(G) > r$.

We have added an arc $e$ to a directed acyclic graph $G$ and have discussed the conditions that the vertices of $G$ adjacent to $e$ should satisfy for the condition $B(G+e)\le B(G)$ to hold. It naturally raises the following question: what conditions should be satisfied in order for this bound to be tight? In a future study  the brushing number of directed graphs with cycles can be discussed in general. More specifically, one could look for the brushing number of tournaments with cycles.
\fi

In Theorem~\ref{th 2} we showed that
if $G$ is a rotational tournament on $n$ vertices then $B(G) = \frac{n^2-1}{8}$.  Rotational tournaments form a subclass of regular tournaments.  We conjecture that the brushing number of regular tournaments can be bounded as follows.
%Based on the observations made on regular tournaments and results obtained about rotational tournaments we present the following conjecture.
\begin{conj}
    If $G$ is regular tournament, then $B(G)\le \frac{n^2-4n+7}{4}$.
\end{conj}

The cleaning model presented in this paper allows any number of brushes to simultaneously traverse an arc of a directed graph.  Imposing edge capacity restrictions would create a new cleaning model for brushing directed graphs; the number of brushes needed to clean a directed graph $G$ under this restricted scenario would serve as an upper bound on our brushing number $B(G)$.

In this paper our focus was
%The cleaning  model that we have considered in this paper is based
on minimising the number of brushes that can be used to clean a given directed graph. It remains as an open problem to explore what is the most efficient cleaning sequence for a directed graph, that is the one that minimises the number of brushes used as well as minimises the time taken to clean the graph when allowing multiple vertices that are ready to fire to do so at the same time step. A way of firing brushes similar to this for undirected graphs is described as \textit{parallel dispersal mode}  in \cite{McKeil2007}.

\section{Acknowledgements}
David Pike acknowledges research support from NSERC Discovery Grant RGPIN-2022-03829.
%\nocite{Sulani2024}

%\bibliographystyle{abbrv}
%\bibliographystyle{elsarticle-num-names}
%\bibliographystyle{elsarticle-num}
\bibliographystyle{plain}
%\bibliography{cas-refs}
\bibliography{paperbib}

\begin{thebibliography}{10}

\bibitem{alspach2004searching}
Brian Alspach.
\newblock Searching and sweeping graphs: a brief survey.
\newblock {\em Le matematiche}, 59(1, 2):5--37, 2004.

\bibitem{BJORNER1991283}
Anders Björner, László Lovász, and Peter~W. Shor.
\newblock Chip-firing games on graphs.
\newblock {\em European Journal of Combinatorics}, 12(4):283--291, 1991.

\bibitem{BonatoNowakowski2011}
Anthony Bonato and Richard~J. Nowakowski.
\newblock {\em The game of cops and robbers on graphs}, volume~61 of {\em
  Student Mathematical Library}.
\newblock American Mathematical Society, Providence, RI, 2011.

\bibitem{bondy2008graph}
J.~A. Bondy and U.~S.~R. Murty.
\newblock {\em Graph theory}, volume 244 of {\em Graduate Texts in
  Mathematics}.
\newblock Springer, New York, 2008.

\bibitem{bryant2014brushing}
Darryn Bryant, Nevena Franceti{\'c}, Przemys{\l}aw Gordinowicz, David~A Pike,
  and Pawe{\l} Pra{\l}at.
\newblock Brushing without capacity restrictions.
\newblock {\em Discrete Applied Mathematics}, 170:33--45, 2014.

\bibitem{DHP2021}
Danny Dyer, Jared Howell, and Brittany Pittman.
\newblock The watchman's walk problem on directed graphs.
\newblock {\em Australas. J. Combin.}, 80:197--216, 2021.

\bibitem{10.1007/978-3-030-83823-2_112}
Alberto Espuny~D{\'i}az, Viresh Patel, and Fabian Stroh.
\newblock Path decompositions of random directed graphs.
\newblock In Jaroslav Ne{\v{s}}et{\v{r}}il, Guillem Perarnau, Juanjo Ru{\'e},
  and Oriol Serra, editors, {\em Extended Abstracts EuroComb 2021}, pages
  702--706, Cham, 2021. Springer International Publishing.

\bibitem{gross2003handbook}
Jonathan~L. Gross and Jay Yellen, editors.
\newblock {\em Handbook of graph theory}.
\newblock Discrete Mathematics and its Applications (Boca Raton). CRC Press,
  Boca Raton, FL, 2004.

\bibitem{HRW1998}
B.~L. Hartnell, D.~F. Rall, and C.~A. Whitehead.
\newblock The watchman's walk problem: an introduction.
\newblock In {\em Proceedings of the {T}wenty-ninth {S}outheastern
  {I}nternational {C}onference on {C}ombinatorics, {G}raph {T}heory and
  {C}omputing ({B}oca {R}aton, {FL}, 1998)}, volume 130, pages 149--155, 1998.

\bibitem{Sulani2024}
Sulani~Dinya Kavirathne.
\newblock {\em Brushing Directed Graphs}.
\newblock M.Sc. Thesis, Memorial University of Newfoundland, 2024.

\bibitem{https://doi.org/10.1112/plms.12328}
Allan Lo, Viresh Patel, Jozef Skokan, and John Talbot.
\newblock Decomposing tournaments into paths.
\newblock {\em Proceedings of the London Mathematical Society},
  121(2):426--461, 2020.

\bibitem{McKeil2007}
Sable~Gwendalyn McKeil.
\newblock {\em Graph cleaning}.
\newblock M.Sc. Thesis, Dalhousie University, 2007.

\bibitem{messinger2008cleaning}
Margaret-Ellen Messinger, Richard~J Nowakowski, and P~Pra{\l}at.
\newblock Cleaning a network with brushes.
\newblock {\em Theoretical Computer Science}, 399(3):191--205, 2008.

\bibitem{MR3391773Penso}
L.~D. Penso, D.~Rautenbach, and A.~Ribeiro~de Almeida.
\newblock Brush your trees!
\newblock {\em Discrete Appl. Math.}, 194:167--170, 2015.

\end{thebibliography}







% \end{thebibliography}
\end{document}